\documentclass[reqno]{amsproc}
\usepackage[colorinlistoftodos,shadow]{todonotes}
\usepackage{amsmath,amscd,amssymb}
\usepackage{cite}
\usepackage{color}
\usepackage{graphicx}
\usepackage{psfrag,epsfig}
\usepackage{epstopdf}
\usepackage{multirow}
\usepackage{graphicx}
\usepackage{rotating}
\usepackage{fullpage}
\usepackage{subfigure}
\usepackage{enumerate}
\usepackage{comment}
\usepackage{lineno}
\usepackage{algorithm2e}
\usepackage[small]{caption}
\usepackage{graphicx}
\usepackage{tabu} 
\usepackage{float}
\usepackage{mathrsfs}
\usepackage{amssymb}
\usepackage{amsmath}
\usepackage{mathtools}
\usepackage[colorlinks=true, linkcolor=blue, citecolor=blue]{hyperref}
\usepackage{siunitx,array,multirow}

\newcommand{\gradient}[1]{\ensuremath{\nabla{#1}}}

\numberwithin{equation}{section}
\newcommand*\patchAmsMathEnvironmentForLineno[1]{%
  \expandafter\let\csname old#1\expandafter\endcsname\csname #1\endcsname
  \expandafter\let\csname oldend#1\expandafter\endcsname\csname end#1\endcsname
  \renewenvironment{#1}%
     {\linenomath\csname old#1\endcsname}%
     {\csname oldend#1\endcsname\endlinenomath}}%
\newcommand*\patchBothAmsMathEnvironmentsForLineno[1]{%
  \patchAmsMathEnvironmentForLineno{#1}%
  \patchAmsMathEnvironmentForLineno{#1*}}%
\AtBeginDocument{%
\patchBothAmsMathEnvironmentsForLineno{equation}%
\patchBothAmsMathEnvironmentsForLineno{align}%
\patchBothAmsMathEnvironmentsForLineno{flalign}%
\patchBothAmsMathEnvironmentsForLineno{alignat}%
\patchBothAmsMathEnvironmentsForLineno{gather}%
\patchBothAmsMathEnvironmentsForLineno{multline}%
}

\def\EQ#1\EN{\begin{equation}#1\end{equation}}
\def\SEQ#1\SEN{\begin{subequations}#1\end{subequations}}
\def\BA#1\EA{\begin{align}#1\end{align}}
\def\LM#1\EM{\begin{linenomath}#1\end{linenomath}}

\newcommand{\pfx}[2]{\dfrac{\partial{#1}}{\partial{#2}}}

\newcommand{\bss}[1]{\boldsymbol{#1}}

\title[Physics-constrained machine learning]{AdjointNet:~Constraining machine learning models with physics-based codes}

\author[S.~Karra et al.]{Satish~Karra$^{1,*}$, Bulbul~Ahmmed$^{1}$, and Maruti~K.~Mudunuru$^{2}$\\
\\
\textit{{$^{1}$\small Computational Earth Science Group (EES-16),
Earth and Environmental Sciences Division, 
Los Alamos National Laboratory, Los Alamos, NM 87545} \\
{$^{2}$\small Watershed \& Ecosystem Science, Pacific Northwest National Laboratory, Richland, WA 99352 }} }
\thanks{$^*$Corresponding author, \href{mailto:satkarra@lanl.gov}{\texttt{satkarra@lanl.gov}}; Last updated: \today; LA-UR-21-28763; PNNL-SA-158825}
\begin{document}
\maketitle
%
%
%
%
\section*{Abstract} 
Physics-informed Machine Learning has recently become attractive for learning physical parameters and features from simulation and observation data.
However, most existing methods do not ensure that the physics, such as balance laws (e.g., mass, momentum, energy conservation), are constrained.
Some recent works (e.g., physics-informed neural networks) softly enforce physics constraints by including partial differential equation (PDE)-based loss functions but need re-discretization of the PDEs using auto-differentiation. 
Training these neural nets on observational data showed that one could solve forward and inverse problems in one shot.
They evaluate the state variables and the parameters in a PDE.
This re-discretization of PDEs is not necessarily an attractive option for domain scientists that work with physics-based codes that have been developed for decades with sophisticated discretization techniques to solve complex process models and advanced equations of state.
This paper proposes a physics constrained machine learning framework, AdjointNet, allowing domain scientists to embed their physics code in neural network training workflows. 
This embedding ensures that physics is constrained everywhere in the domain. 
Additionally, the mathematical properties such as consistency, stability, and convergence vital to the numerical solution of a PDE are still satisfied.
We show that the proposed AdjointNet framework can be used for parameter estimation (and uncertainty quantification by extension) and experimental design using active learning.
The applicability of our framework is demonstrated for four cases -- (1) flow in a homogeneous porous medium, (2) data assimilation for homogeneous porous media flow, 
(3) flow in a heterogeneous porous medium, and (4) cavity flow using the Navier-Stokes equation. 
Results show that AdjointNet-based inversion can estimate process model parameters with reasonable accuracy.
These examples demonstrate the applicability of using existing software with no changes in source code to perform accurate and reliable inversion of model parameters. 
\\ \\
\noindent\textbf{Keywords:}~machine learning, physics constraints, partial differential equations, adjoint, neural networks, inversion. 
%
\section{Introduction}
\label{sec:introduction}
State-of-the-art physics-informed machine learning (PIML) \cite{harp2021feasibility,kashinath2021physics,yan2021physics} or knowledge-guided machine learning (KGML) \cite{karpatne2017theory,hanson2020predicting,khandelwal2020physics} approaches have been primarily data-driven. 
PIML model development involves data generation from physics-based codes.
Then, unsupervised (e.g., NMFk, NTFk, PCA, autoencoders) \cite{vesselinov2019unsupervised} or supervised (e.g., deep neural networks) \cite{harp2021feasibility} learning methods are used to train on
this data to obtain state variables or quantities of interest. 
These approaches decompose complex coupled processes into various features/signals that help domain scientists identify the unique physical mechanisms \cite{cichocki2009nonnegative,vesselinov2019unsupervised}. 
Furthermore, popular PIML methods such as physics-informed neural networks (PINNs) \cite{raissi2019physics,tartakovsky2020physics,lu2021deepxde} softly (e.g., in $L^2$-norm) ensure physics constraints by training the PDE variables as neural networks.
Auto-differentiation (e.g., in \textsf{TensorFlow} \cite{abadi2016tensorflow}, \textsf{PyTorch} \cite{paszke2017automatic,paszke2019pytorch}) is used to discretize the PDE. 
One then uses a weighted loss function (e.g., based on mean squared error) that linearly combines the loss due to PDE residual and discrepancies from
data.
This loss function ensures that the governing laws are `reasonably satisfied' although not entirely accurate locally at mesh/collocation points in the domain \cite{jagtap2020conservative}.

Existing PIML and KGML approaches also require large volumes of data to make reasonably accurate predictions. 
As a result, ML model development becomes computationally intensive and time-consuming as large data needs to be generated. 
For example, it takes several hours to a day to run a single at-scale subsurface reservoir simulation \cite{mudunuru2020physics} with
degrees-of-freedom in the $\mathcal{O}(10^7)$ on state-of-the-art high-performance computing (HPC) machines \cite{NERSC2021,OLCF2021,ALCF2021}. 
For this reason, existing PIML/KGML methods are not ideal for approximating complex and natural systems, which require at least hundreds of realizations as training data.
Furthermore, for applications where heterogeneity is present in the systems, for instance, different materials in different parts of the domain, optimizing the number of parameters can be significant. 
For instance, in subsurface flow applications, properties such as permeability, a measure of how much a porous material can allow flow, can be different in each grid cell of the model. 
In such scenarios, traditional ML or PIML methods require Monte Carlo approaches (e.g., MCMC) to sample from the parameter spaces and generate data \cite{santos2020modeling}. 
These requirements make data generation even more challenging.

\begin{figure}[htbp!]
  \centering
    \subfigure[State-of-the-art ML workflow]
    {\includegraphics[width = 0.9\textwidth]
    {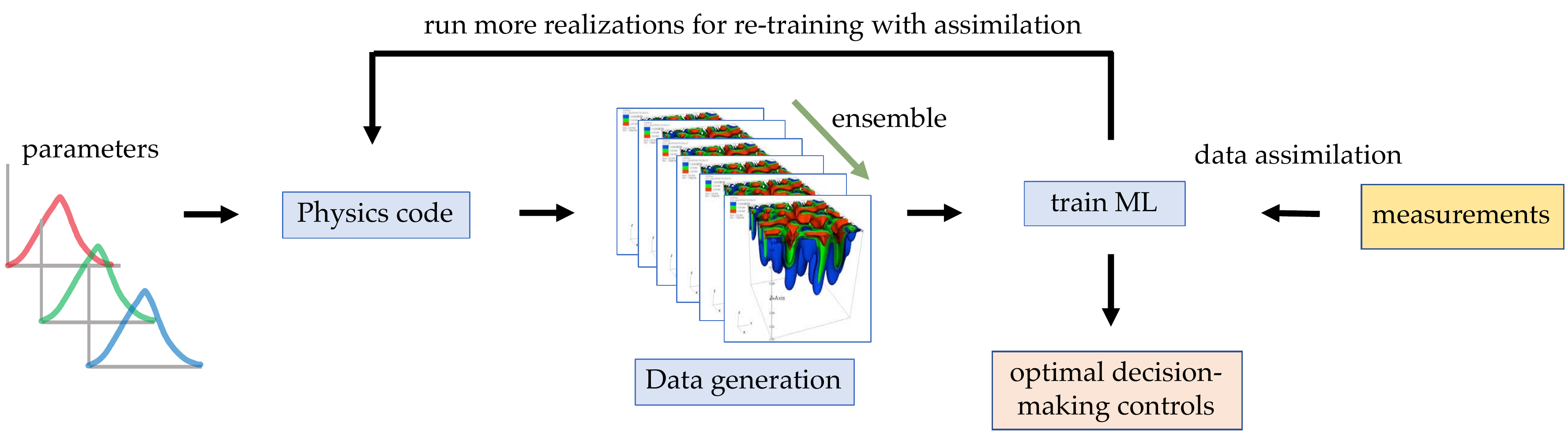}}
  %
  %
    \subfigure[Proposed workflow]
    {\includegraphics[width = 0.6\textwidth]
    {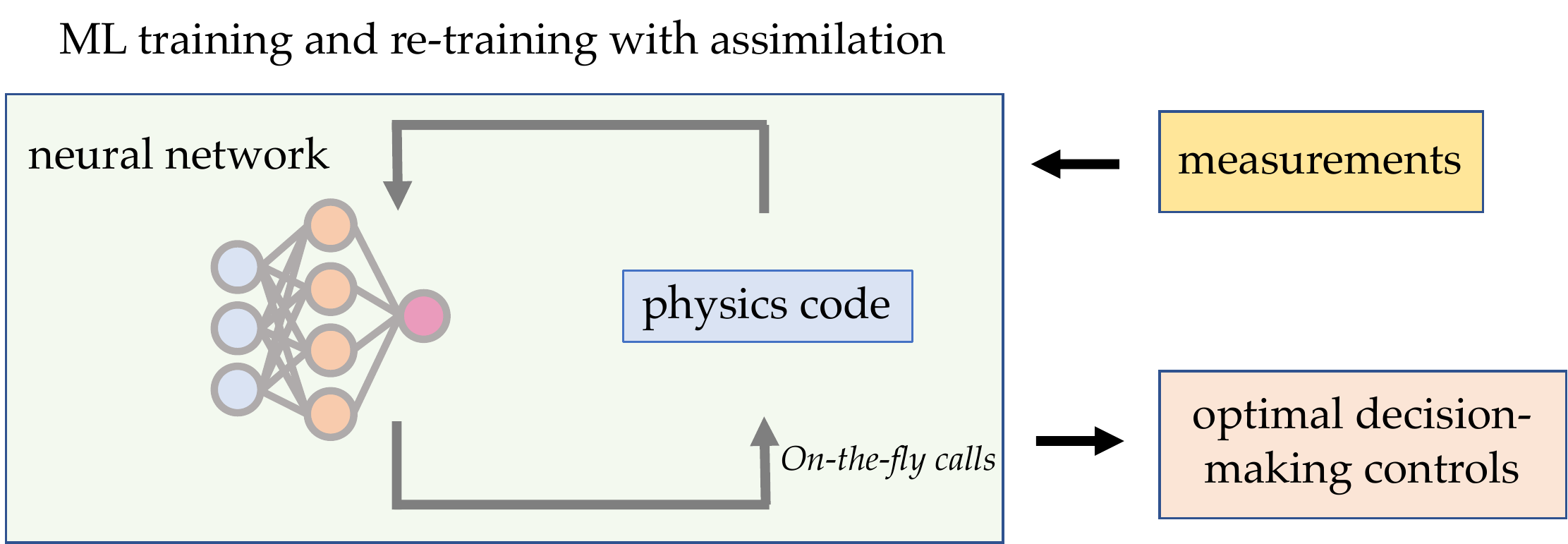}}
  \caption{\textbf{AdjointNet framework}:~Comparison between the state-of-the-art ML workflow with the proposed workflow.
  A major advantage is that the proposed workflow requires minimal simulations, as it calls the physics-based code on the fly, to perform data assimilation and ML training.
  \label{Fig:workflow}}
\end{figure}

Most PIML/KGML methods do not include simulation data and experimental observations in the loop during the ML training process. 
Training is usually a two-step process, where an ML model is first built from simulations. 
Then real-time data streams are integrated later to infer system decision parameters (as shown in Fig.~\eqref{Fig:workflow}(a)).
If new data is assimilated beyond the training range, one needs to generate more training data from simulations and then re-train the ML model. 
Our proposed method overcomes this serious technical gap by directly incorporating physics-based codes into the ML training procedure.
We can use our approach to perform ML while ensuring that a given physics constraint (e.g., the balance of mass, the balance of momentum) is satisfied everywhere, as modeled by a physics-based code. 
\textit{Our innovation is the use of adjoint sensitivities, which can be obtained directly from physics calculations, to calculate the gradients of the loss function in the training process of a neural network.}
Since one needs to solve the forward problem using the physics code to calculate these adjoint sensitivities, the underlying physics will be `automatically' satisfied and constrained. 
Due to the `correct' constraining of physics in the entire domain, one can still get highly predictive ML models with limited, sparsely-sampled, and noisy data.

Our computational approach (see Fig.~\eqref{Fig:workflow}(b)), AdjointNet, allows domain scientists to incorporate simulators into scalable ML workflows without any loss of physics (up to the accuracy of a simulator) and simultaneously integrate observational data.
Moreover, our method can directly assimilate real-time data streams in the training step. 
The trained ML models can then be used to forecast system behavior or inform system optimization.
Furthermore, this approach has the advantage that one can invert the unknown parameters in training the ML model to measured data.
We can use these parameters to forecast the system behavior towards further decision-making. 
In this paper, we demonstrate the utility of this AdjointNet framework by applying it to inversion and learn model parameters in porous media flow and Navier-Stokes flow problems, where the balance laws -- the balance of mass, the balance of momentum -- are solved. The following specific examples are solved: (1) flow in a homogeneous porous medium, (2) data assimilation for homogeneous porous media flow, 
(3) flow in a heterogeneous porous medium, and (4) cavity flow using the Navier-Stokes equation.

Below, we summarize the novelty of the proposed AdjointNet framework compared to current state-of-the-art:
\begin{enumerate}
 \item \textbf{Incorporating existing codes and extendability}:~Our framework can incorporate any existing physics-based codes such as \texttt{PFLOTRAN}\cite{hammond2014evaluating,lichtner2015pflotran}, \texttt{E3SM} \cite{golaz2019doe}, \texttt{OpenFOAM} \cite{jasak2007openfoam}, \texttt{SWAT} \cite{arnold2012swat}, \texttt{PRMS}\cite{leavesley1995precipitation}, \texttt{HOSS} \cite{knight2020hoss}, and \texttt{WRF-Hydro} \cite{sampson2018wrf}.
 For instance, these codes can simulate complex coupled processes such as flow, chemical and thermal transport, and mechanics for various applications such as climate, subsurface energy, and carbon and nuclear waste storage, etc. 
 Most of the above codes have undergone rigorous verification and validation over the years.
 As a result, we can simulate, estimate, and invert for material or conceptual properties with no changes in existing codes reliably.
 This provides a significant advantage over existing methods such as PINNs \cite{raissi2019physics}, where process models and equations of state need to be re-written to perform the inversion.
 \item \textbf{Balance laws}:~Even though the trained PIML models seemingly give reasonable predictions on quantities of interest such as pressures, displacements, and temperatures, there are no guarantees that the governing laws of physics (e.g., the balance of mass, energy, or momentum) will truly be satisfied by the machine learning (ML) model.
 That is, the model can violate these laws locally or globally, or both, as reported by the developers of PINNs \cite{jagtap2020conservative}. 
 For example, a significant limitation of PINNs \cite{fuks2020limitations,wang2020understanding,jagtap2020conservative} is the accuracy of the solution and associated local balance.
 In PINNs, the absolute error does not go below $\mathcal{O}(10^{-5})$ \cite{jagtap2020conservative}.
 This is because of the reduced accuracy in solving the high-dimensional non-convex optimization problem, resulting in local minima.
 Since the AdjointNet framework directly calls the existing codes, which ensure that the balance laws are satisfied up to machine precision, the final trained model will always be locally conservative \cite{eymard2000finite,hughes2000continuous}. 
 \item \textbf{Scalability}:~Over the decades, some of these codes mentioned above have been parallelized to scale on leadership-class HPC machines \cite{hammond2014evaluating,lichtner2015pflotran,golaz2019doe}.
 Our AdjointNet framework can be scalable as it uses the existing HPC codes under the hood.
 Moreover, the proposed algorithm uses discrete adjoint-based methods, which are also scalable \cite{towara2015mpi,zhang2019petsc}.
 Hence, the overall approach can be scaled on the next generation HPC machines geared towards exascale computing (e.g., Perlmutter, Aurora) \cite{messina2017exascale,stevens2019aurora,yang2020accelerate}.
 This provides a significant advantage over other PIML methods (e.g., PINNs).
 \item \textbf{Agnostic to underlying numerical discretization}:~AdjointNet works with any numerical method implemented in a physics code.  
 Here, we demonstrate the utility with two physics codes implemented with different numerical methods--finite volume and finite difference--but can be extended to any numerical method of choice.
 \item \textbf{Differentiable programming}:~Recently, programming languages like Julia have made differentiable programming more accessible \cite{innes2019differentiable}.
 However, the governing equations in physics codes have to be re-written in Julia to take advantage of this capability.
 As mentioned previously, this re-writing of codes may not be attractive for domain scientists. 
 AdjointNet leverages the infrastructure in \textsf{Tensorflow} or \textsf{PyTorch} to perform auto-differentiation on the neural networks but does not require differentiable programming of the underlying physics codes.
 \item \textbf{Computational cost}:~Another major limitation of PINNs in addition to inaccuracy is the significant training cost associated with deep neural networks and long-time integration of the PDEs \cite{jagtap2020conservative}.
 This is because the addition of a PDE-based loss function results in a stiff and non-convex optimization problem \cite{xu2020inverse,xu2020physics}.
 Our AdjointNet framework does not do this but poses ML training as a PDE-constrained optimization problem.
 Furthermore, since we can directly integrate existing HPC physics-based codes, we can obtain the final trained model faster, thus allowing us to work with more realistic and large-scale problems.
 \item \textbf{Desirable mathematical properties}:~Most of the existing physics codes satisfy the desirable properties such as structure-preserving, coercivity, error estimates, consistency, convergence, accuracy, and stability \cite{zienkiewicz2005finite,eymard2000finite}.
 Whichever of these properties are satisfied by the underlying PDE solver used in the AdjointNet framework, the final trained model will also satisfy that.
 On the other hand, existing PIML methods do not ensure this.
 \item \textbf{eXplainable Artificial Intelligence (XAI)}:~Typically, numerical methods like finite element, finite volume, etc., \cite{zienkiewicz2005finite,eymard2000finite} are employed to obtain the solution of a PDE.
 These numerical methods are explainable, which results in the trustworthiness of the obtained solution.
 On the other hand, existing PIML methods such as PINNs use neural networks to compute the solution of the variables in a PDE.
 These neural networks are difficult to interpret, and XAI methods \cite{dovsilovic2018explainable,adadi2018peeking} (e.g., SHAPley values \cite{sundararajan2020many}, LRP \cite{montavon2019layer}, deep Taylor decomposition \cite{montavon2017explaining}, LIME \cite{ribeiro2016model}) are needed to solve the solution.
 When neural networks are employed to obtain the solution, this interpretation process needs to be performed for each scenario of interest (e.g., changes in boundary conditions, initial conditions, or domain).
 However, note that most of the PIML related works have rarely used XAI methods to explain why such a method works.
 \item \textbf{Reproducability}:~Typically, in any PIML method, the solution to the PDE depends on the trainable weights that need to be initialized randomly \cite{kumar2017weight,nagarajan2019generalization,zhang2020type}.
 As a result, the procedure of weights initialization plays a predominant role in the obtained solution, leading to inconsistency and hence may not be reproducible.
 However, our AdjointNet is not affected by such initialization as we rely on the existing physics codes to always provide reproducible solutions of the PDE. 
\end{enumerate}

The outline of our paper is as follows:~Sec.~\eqref{sec:formulation} presents the mathematical formulation and the linkage between physics-based codes, which is the AdjointNet framework. 
Section~\eqref{sec:examples} presents examples of our proposed approach.
Conclusions are drawn in Sec.~\eqref{sec:conclusions}.

\section{AdjointNet Framework}\label{sec:formulation}

Machine learning methods (e.g., neural networks) typically involve minimizing a loss function $\mathcal{L}$. 
We add a constraint that the PDE of interest needs to be satisfied in our underlying formulation. 
This turns into a PDE-constrained optimization problem, that is, 
\begin{align}  
  \textrm{min}~ \mathcal{L} \label{eq:loss} \\
  \textrm{s.t.} ~ \mathcal{F}\left(\boldsymbol{u}(\boldsymbol{x}, t); \boldsymbol{p}\right) = \boldsymbol{0} \label{eq:pde}
\end{align}
where $\boldsymbol{x}$ is the position, $\boldsymbol{u}$ is the variable of interest (e.g., pressure, temperature, concentration), $\boldsymbol{p}$ is the vector of parameters (e.g., permeability, thermal conductivity, diffusivity), and $t$ is time.
Here,  Eq.~\eqref{eq:pde} is the discrete form of the PDE that is solved by a physics-based code. For instance, in the case of porous media, the parallel code \texttt{PFLOTRAN} solves the discrete form of the non-linear diffusion PDE \cite{lichtner2015pflotran}:
\BA
 \pfx{\left[\phi\rho(u) u\right]}{t} - \mathrm{div}\left[\frac{\rho(u) k}{\mu} \gradient{u}\right] = 0 
\EA
where $u$ is the fluid pressure, $k$ is the permeability parameter, $\phi$ is porosity, $\mu$ is fluid viscosity, and $\rho$ is the fluid density. 
PDE discretization in \texttt{PFLOTRAN} is performed using a two-point flux finite volume method with backward Euler for time-stepping. 
The resulting nonlinear algebraic equations are solved using a Newton-Krylov solver.

If we assume that the parameters $\boldsymbol{p}$ are to be learnt from a set of observational data $\boldsymbol{u}_{\mathrm{obs}}$, then one can use a loss function of $\mathcal{L} = \|u-u_{obs}\|$. 
Also, suppose that $\boldsymbol{p}$ is a ML model such as a neural network, i.e., $\bss{p} = NN(\bss{x})$. 
If $W^{(k)}$ and $b^{(k)}$ are the weights and biases for this neural network $NN(\boldsymbol{x})$, then the gradients $\dfrac{\partial \mathcal{L}}{\partial W^{(k)}}$, $\dfrac{\partial \mathcal{L}}{\partial b^{(k)}}$ are needed.
These gradients must be supplied to the optimization algorithm that minimizes the loss function; for instance, stochastic gradient descent is a commonly used algorithm.  
Now, using the chain rule, one can write these gradients as
\begin{align}
  \dfrac{\partial \mathcal{L}}{\partial W^{(k)}} = \dfrac{\partial \mathcal{L}}{\partial \boldsymbol{u}} \dfrac{\partial \boldsymbol{u}}{\partial \boldsymbol{p}} \dfrac{\partial \boldsymbol{p}}{\partial W^{(k)}}, 
  ~\dfrac{\partial \mathcal{L}}{\partial b^{(k)}} = \dfrac{\partial \mathcal{L}}{\partial \boldsymbol{u}} \dfrac{\partial \boldsymbol{u}}{\partial \boldsymbol{p}} \dfrac{\partial \boldsymbol{p}}{\partial b^{(k)}}	\label{eq:lossfunctionbias}
\end{align}
Notice that in both sets of gradients, we need the sensitivity $\dfrac{\partial \boldsymbol{u}}{\partial \boldsymbol{p}}$, calculated through the adjoint form of Eq.~\eqref{eq:pde},  which is obtained by taking the derivative of Eq.~\eqref{eq:pde} with respect to $\bss{p}$:
\BA
\dfrac{\partial \mathcal{F}}{\partial \boldsymbol{p}} + \dfrac{\partial \mathcal{F}}{\partial \boldsymbol{u}}\dfrac{\partial \boldsymbol{u}}{\partial \boldsymbol{p}} = \boldsymbol{0} \label{eq:adjoint}
\EA
We have dropped time $t$ in Eq.~\eqref{eq:adjoint} assuming steady-state, for the sake of illustration.
However, in the case of transient problems, if one discretizes the PDE with time first, we get a similar form for the adjoint equation Eq.~\eqref{eq:adjoint}.
In each epoch of the optimization loop, one needs to solve Eq.~\eqref{eq:adjoint} to get the sensitivity. The other terms in the chain rule, the derivatives $\dfrac{\partial \mathcal{L}}{\partial \boldsymbol{u}}$, $\dfrac{\partial \boldsymbol{p}}{\partial b^{(k)}}$, $\dfrac{\partial \boldsymbol{p}}{\partial W^{(k)}}$ can be obtained directly through auto-differentiation; ML packages like \textsf{TensorFlow} or \textsf{PyTorch} provide them.
One can use the discrete adjoint method \cite{giles2003algorithm} to solve Eq.~\eqref{eq:adjoint}. 

Discrete adjoint methods involve solving the forward problem and checkpointing the solution; then the adjoint (backward) solve is performed by integrating Eq.~\eqref{eq:adjoint} using the saved forward solution at checkpoints. 
Hence, the cost of an adjoint solve is approximately the cost of a forward solve, which is independent of the number of parameters. 
For this reason, discrete adjoint methods are popular in the area of design, for example, in shape optimization of airfoils in the aerospace industry \cite{nadarajah2007optimum,rumpfkeil2010optimal}. 
Since in the discrete adjoint method to solve Eq.~\eqref{eq:adjoint}, one must solve the forward problem Eq.~\eqref{eq:pde}, we ensure that the physics is constrained within each epoch. 
Alternatively, one can obtain the sensitivities $\dfrac{\partial \boldsymbol{u}}{\partial \boldsymbol{p}}$, by perturbing each parameter in $\bss{p}$, re-solving the PDE Eq.~\eqref{eq:pde}, and then calculating the gradients numerically. 
Even for this approach, since one needs to solve the PDE, the physics is thus constrained in each epoch.
However, the cost of evaluating gradients with the proposed approach is $N_p +1$ forward solves, where $N_p$ is the number of parameters. Thus as the unknown parameters increase, this approach becomes relatively expensive.
An approach to overcome this problem is to integrate prior information (e.g., domain expertise, local and global sensitivity analysis, mutual information theory) with unsupervised learning (e.g., $k$-mean clustering) to reduce the dimensions of unknown parameters before applying our AdjointNet framework.
This provides an efficient way to parameterize the high-dimensional space of unknown parameters before performing inversion efficiently.

\section{Numerical Examples}
\label{sec:examples}
In this section, we present representative numerical examples to demonstrate the applicability of our AdjointNet framework for parameter inversion and data assimilation.
\subsection{Porous media flow:~Inversion of homogeneous permeability under sparse data}
\label{subsec:homogeneous_case}
Through a synthetic case, we perform porous media simulations using \texttt{PFLOTRAN} to invert for homogeneous permeability.
\texttt{PFLOTRAN} is an open-source, state-of-the-art massively parallel subsurface flow and reactive transport code \cite{hammond2014evaluating,lichtner2015pflotran}.
We use the pressure data sampled at different locations in the domain as ground truth to perform the inversion.
The model domain is $100 \times 1 \times 1$\,\si{\cubic\metre} in dimension and consists of one layer. 
A uniform mesh of 100 grid cells is used to generate data. 
The permeability value $10^{-14}$~m$^2$ is assigned for all grids, which is assumed to be the ground truth.
The initial condition for the \texttt{PFLOTRAN} simulation includes the pressure of 1.5$\times 10^5$ Pa throughout the model domain.
Flow is driven by pressure boundary conditions,  1.0$\times 10^6$ Pa and 1.5$\times 10^5$ Pa at the left and right faces, respectively.
\begin{figure}
  \centering
    \subfigure[Permeability estimation vs. epochs]
    {\includegraphics[width = 0.49\textwidth]
    {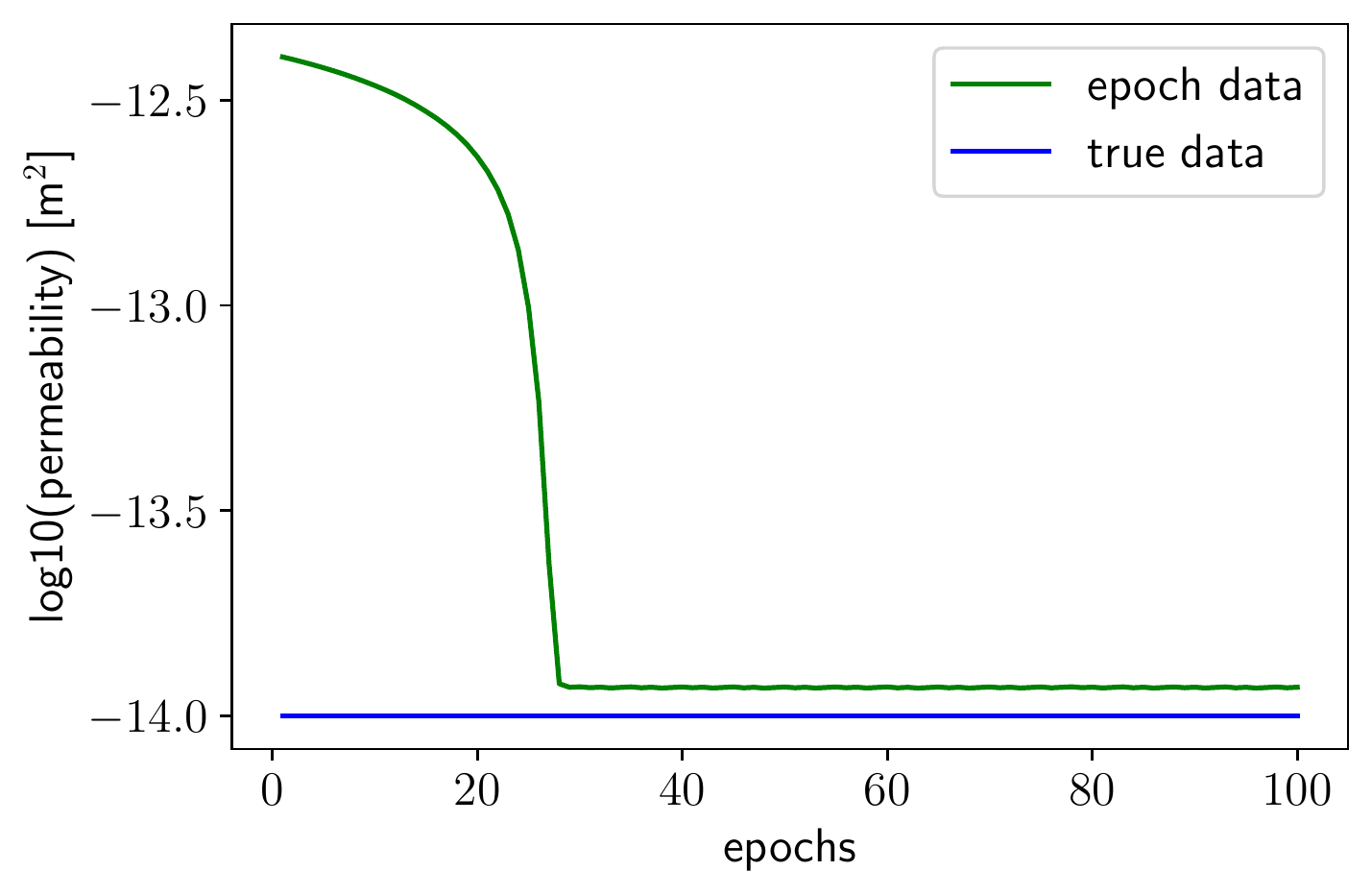}}
  %
  %
    \subfigure[Pressure profile prediction and comparison]
    {\includegraphics[width = 0.49\textwidth]
    {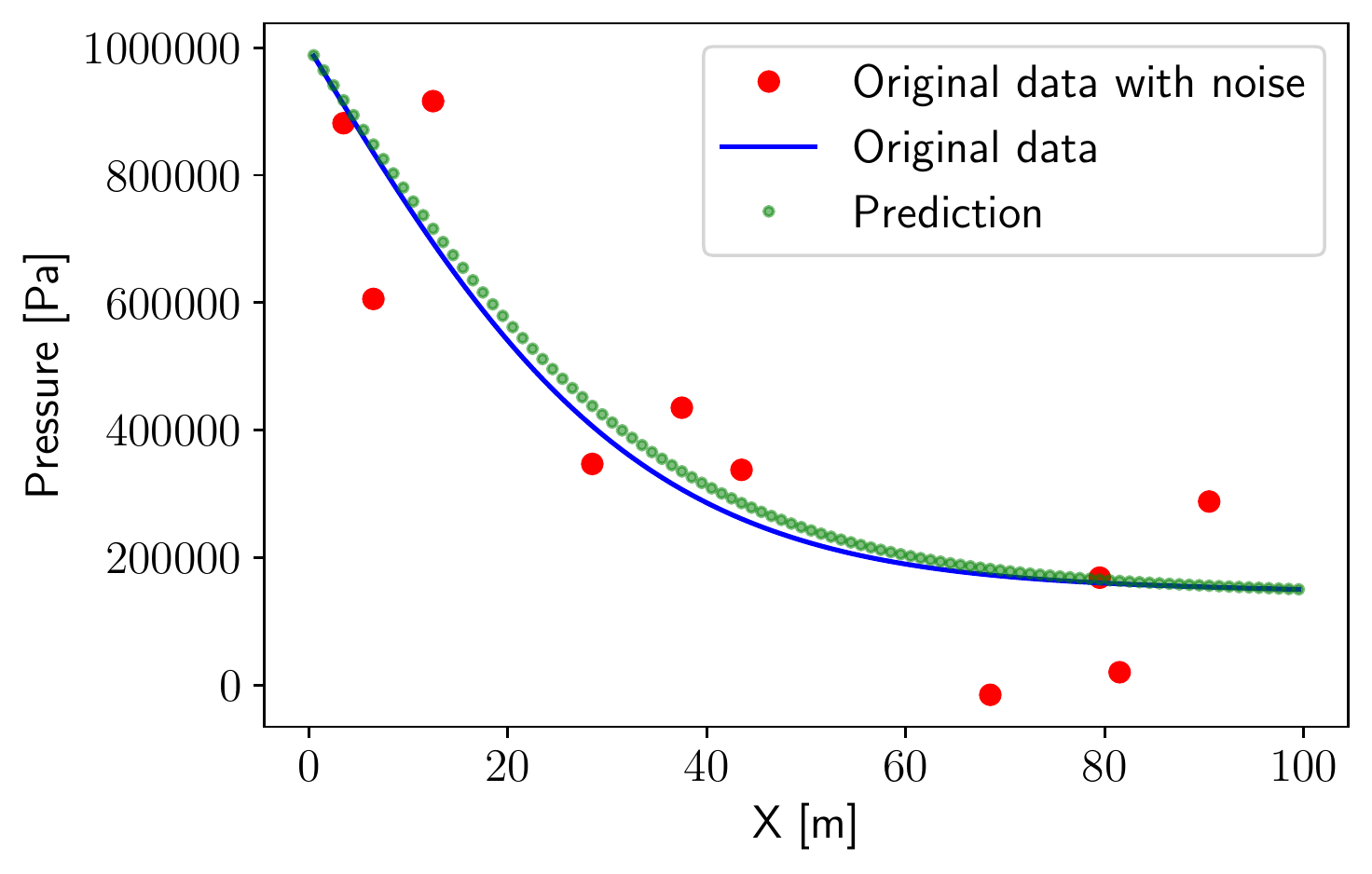}}
  \caption{\textbf{AdjointNet for inversion of homogeneous permeability:}~Training of permeability using AdjointNet. 
  As the model trains to the pressure data (with noise added) as shown in (b), the model learns the permeability as shown in (a). (b) also compares the prediction using the trained permeability against the original pressure data.
  \label{Fig:perm_inversion}}
\end{figure}

To check the validity of the proposed framework, first, we generate ground truth pressure data using the homogeneous numerical simulation.
Next, we sparsely sampled ten out of 100 data points for training the permeability neural network. 
Random noise (with a magnitude of $\approx 0.1$ MPa) is added to the training data to perform inversion under noisy data. 
A permeability neural network is created that has two layers containing 256 neurons in each layer. 
The loss function is formulated based on the mean-squared error between the true and predicted pressure values at these ten points.
The learning rate is set to $10^{-3}$, and the neural network weights and biases are initialized based on Glorot uniform initialization.
The homogeneous permeability neural network is trained for 100 epochs, and convergence is reached within 30 epochs as shown in Fig.~\eqref{Fig:perm_inversion}(a).
During each epoch, the neural network is updated as the flow process is simulated by \texttt{PFLOTRAN}.
Neural network provides gradient of permeability with respect to weights $\dfrac{\partial \bss{p}}{\partial W}$ while the numerical model provides $\dfrac{\partial \mathcal{L}}{\partial \boldsymbol{u}}$ and $\dfrac{\partial \boldsymbol{u}}{\partial \boldsymbol{p}}$.
The loss value plateaued to 0.0015.
Figure~\eqref{Fig:perm_inversion}(a) shows that the AdjointNet estimated permeability is very close to the ground truth.
Figure~\eqref{Fig:perm_inversion}(b) compares the predicted pressures to ground truth and shows that pressure data predictions are robust to noise compared to actual values in the entire domain.
This case study instills confidence that AdjointNet methodology can provide reasonably accurate permeability estimates even under sparsely sampled data and its model predictions are robust to observational noise.

\begin{figure}
  \centering
    \subfigure[Phase~1]
    {\includegraphics[width = 0.49\textwidth]
    {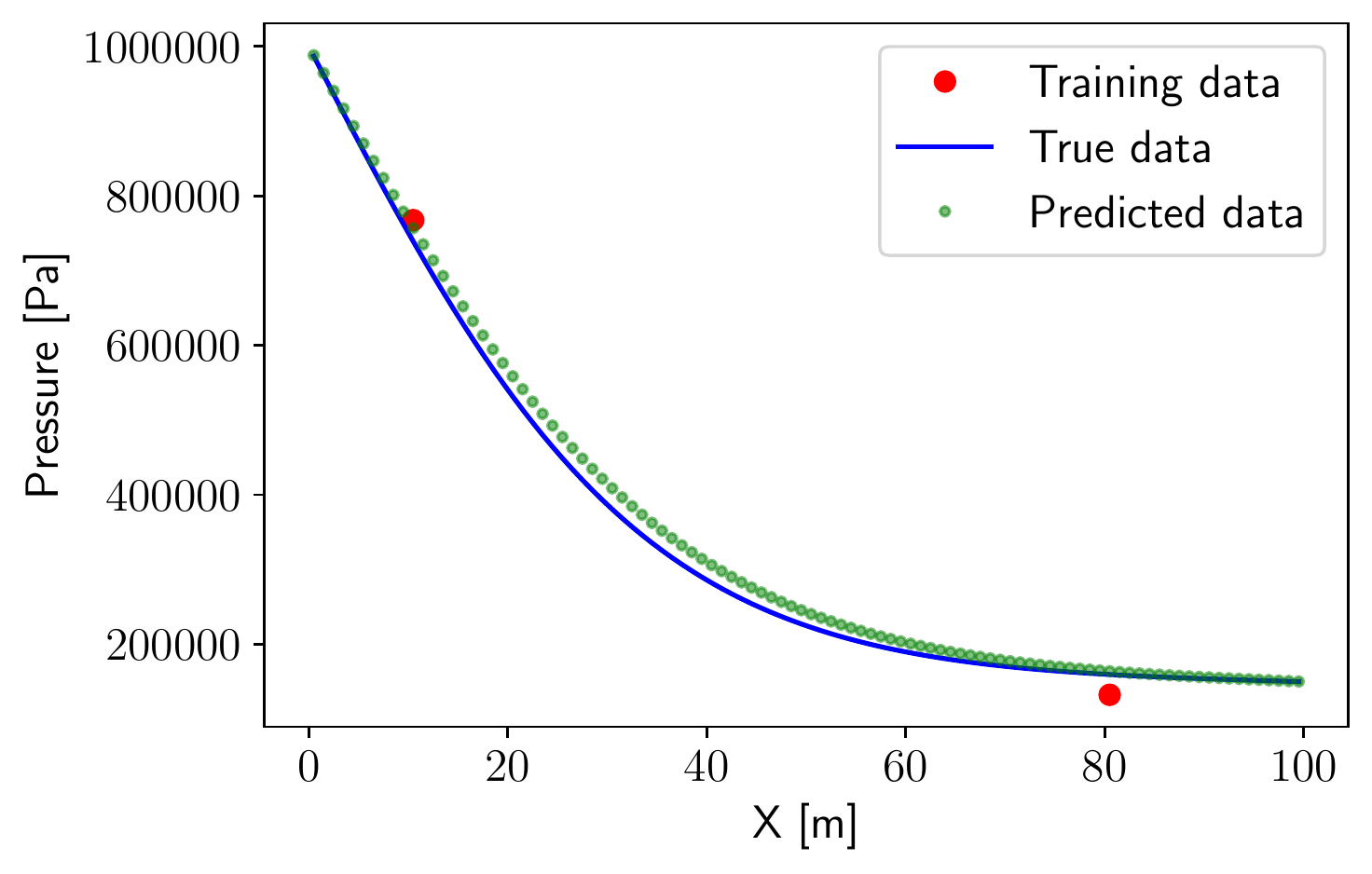}}
  \subfigure[Phase~2]
  {\includegraphics[width = 0.49\textwidth]
  {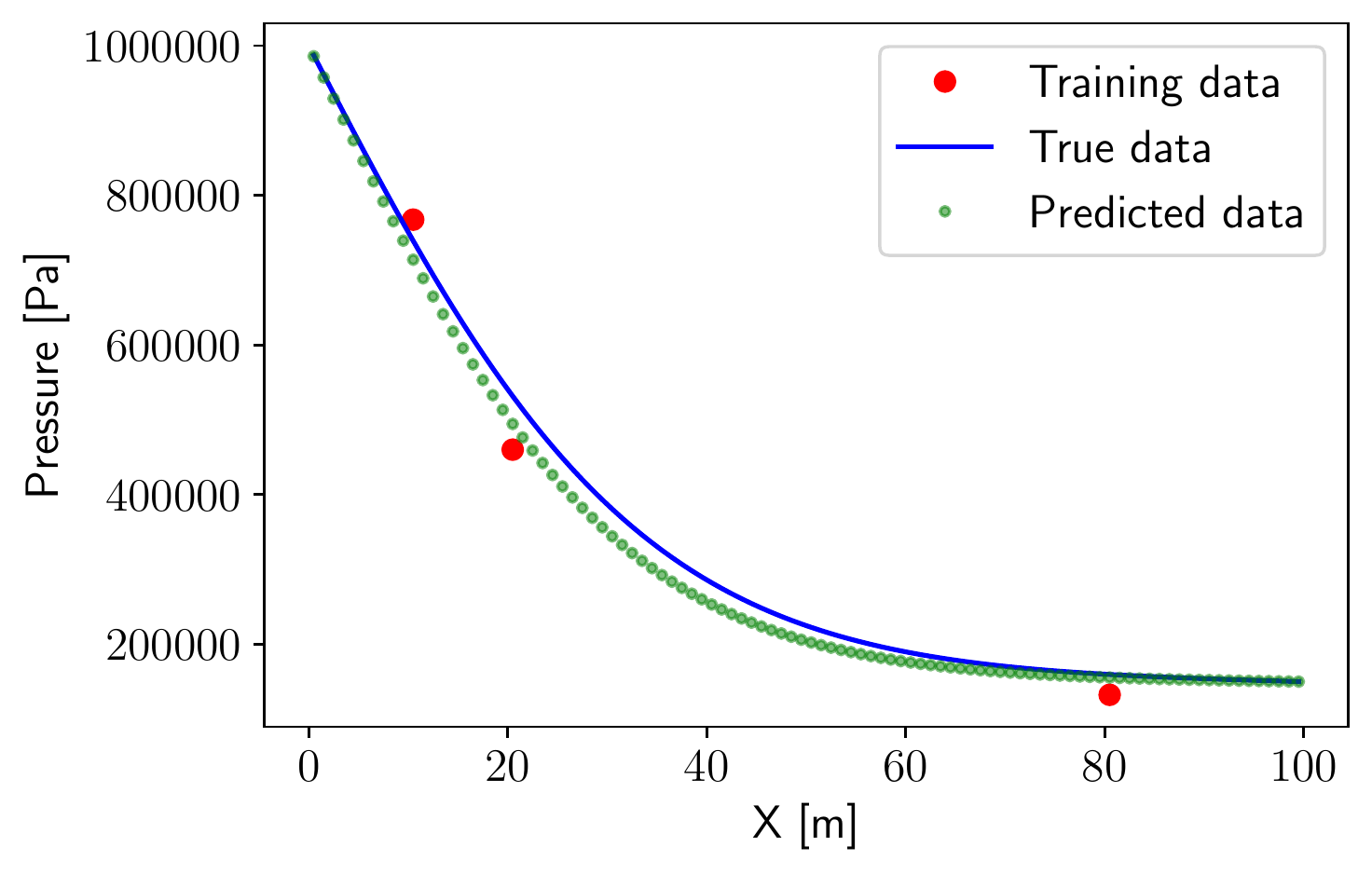}}
  \vspace{0.1in}
  \subfigure[Phase~3]
  {\includegraphics[width = 0.49\textwidth]
  {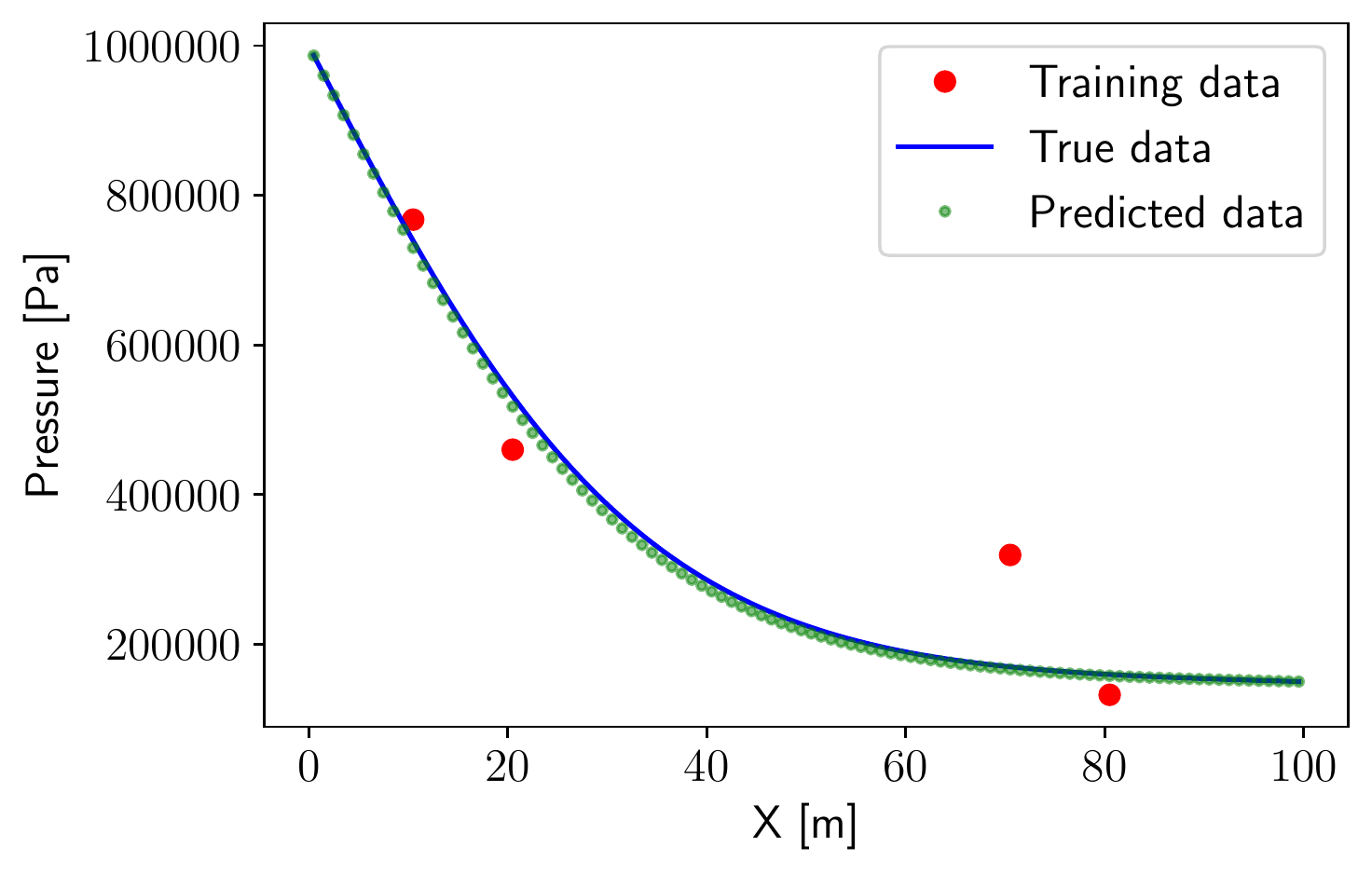}}
  \subfigure[Permeability training]
  {\includegraphics[width = 0.49\textwidth]
  {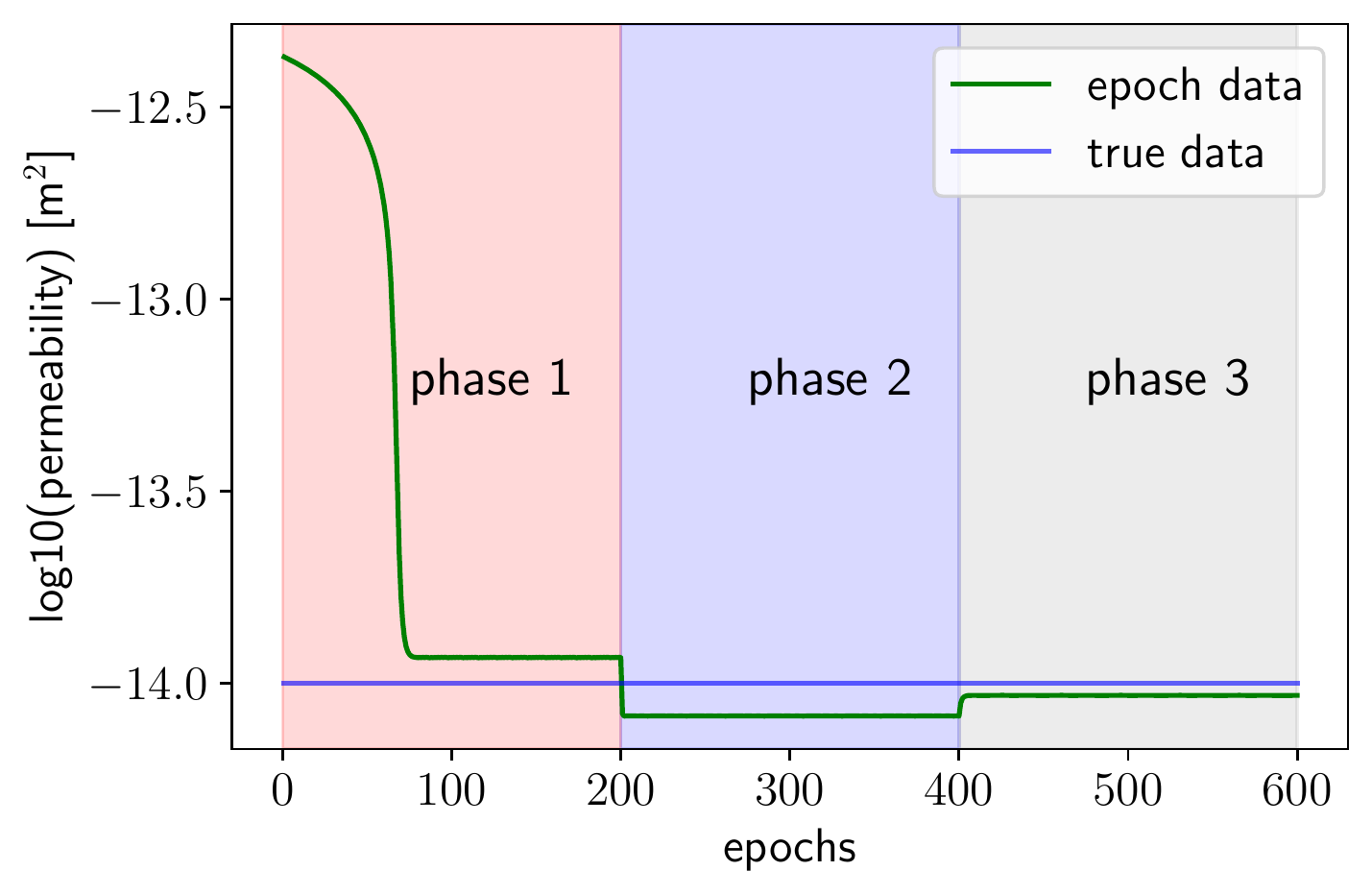}}
  \caption{\textbf{AdjointNet for data assimilation with homogeneous permeability:}~Training of permeability using a neural network as data is assimilated in three sequential phases shown in (a), (b), and (c). The trained permeability is shown in (d). 
  \label{Fig:perm_inversion2}}
\end{figure}

\subsection{Porous media flow:~Data assimilation}
This case study shows a proof-of-concept towards data assimilation applications.
The simulation setup is the same as the previous use case. However, the homogeneous permeability neural network is trained in three sequential phases (Phase~1, Phase~2, and Phase~3). 
The neural network is first trained against two randomly selected pressure data points with added random noise (with a magnitude of $0.1$ MPa) in Phase~1 (Fig.~\eqref{Fig:perm_inversion2}(a)). 
This is followed by further training by assimilating an additional data point in each of Phase~2 (Fig.~\eqref{Fig:perm_inversion2}(b)) and Phase~3 (Fig.~\eqref{Fig:perm_inversion2}(c)). 
Training in each phase is run for 200 epochs, which totals 600 epochs.
Figure~\eqref{Fig:perm_inversion2}(d) shows the permeability predictions for three phases and compares them with the ground truth.
It is evident that as we progressively add more pressure data, we can better constrain permeability.
As a result, Phase 3 permeability prediction is closer to the ground truth than Phase 2 and Phase 1 (see Fig.~\eqref{Fig:perm_inversion2}(d)).

\subsection{Porous media flow: Heterogeneous permeabilities}
\begin{figure}
  \centering
    {\includegraphics[width = 0.75\textwidth]
    {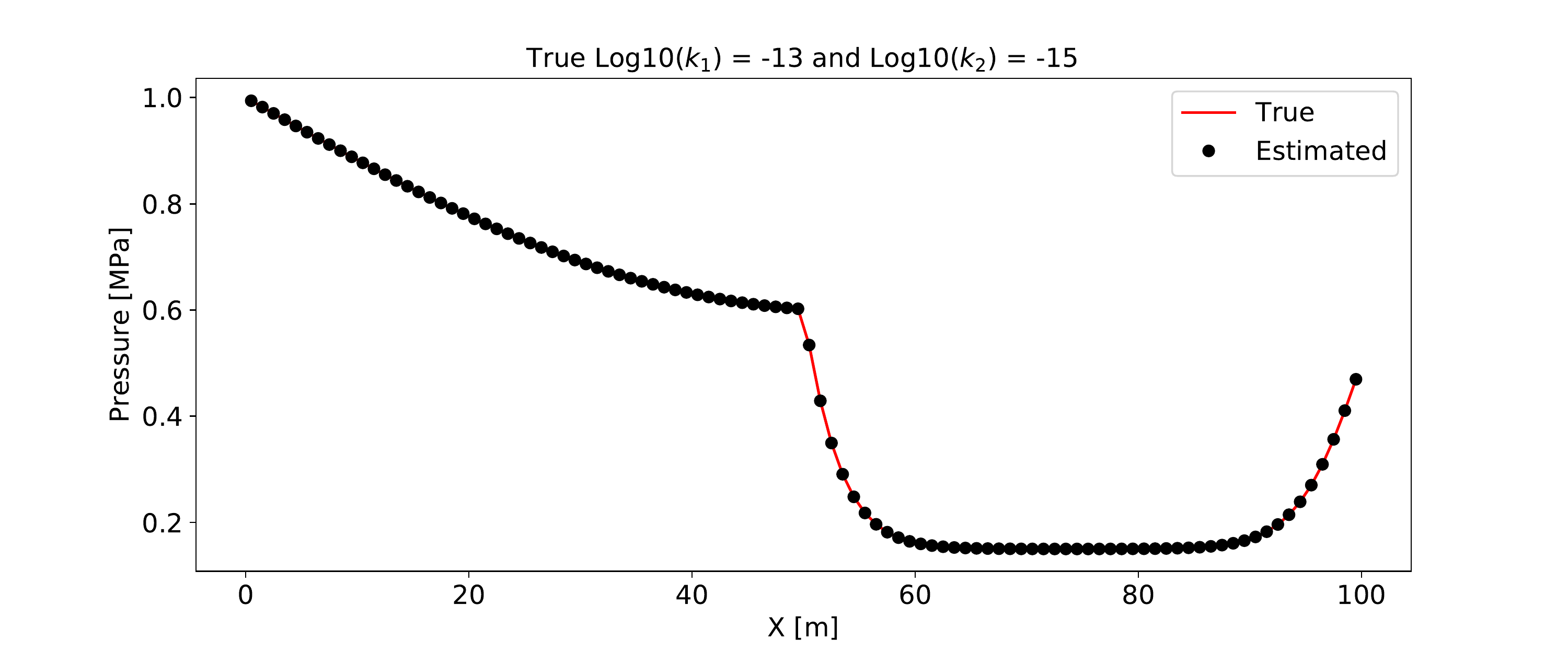}}
  \caption{\textbf{AdjointNet for inversion of heterogeneous permeability:}~This figure shows the actual and AdjointNet predicted pressure distribution over the model domain.
  AdjointNet can capture the pressure profiles with reasonably good accuracy even under discontinuities.}
  \label{fig:hetero-prediction-vs-real}
\end{figure}

\begin{figure}
  \centering
    \subfigure[$k_1$ training]
    {\includegraphics[width = 0.49\textwidth]
    {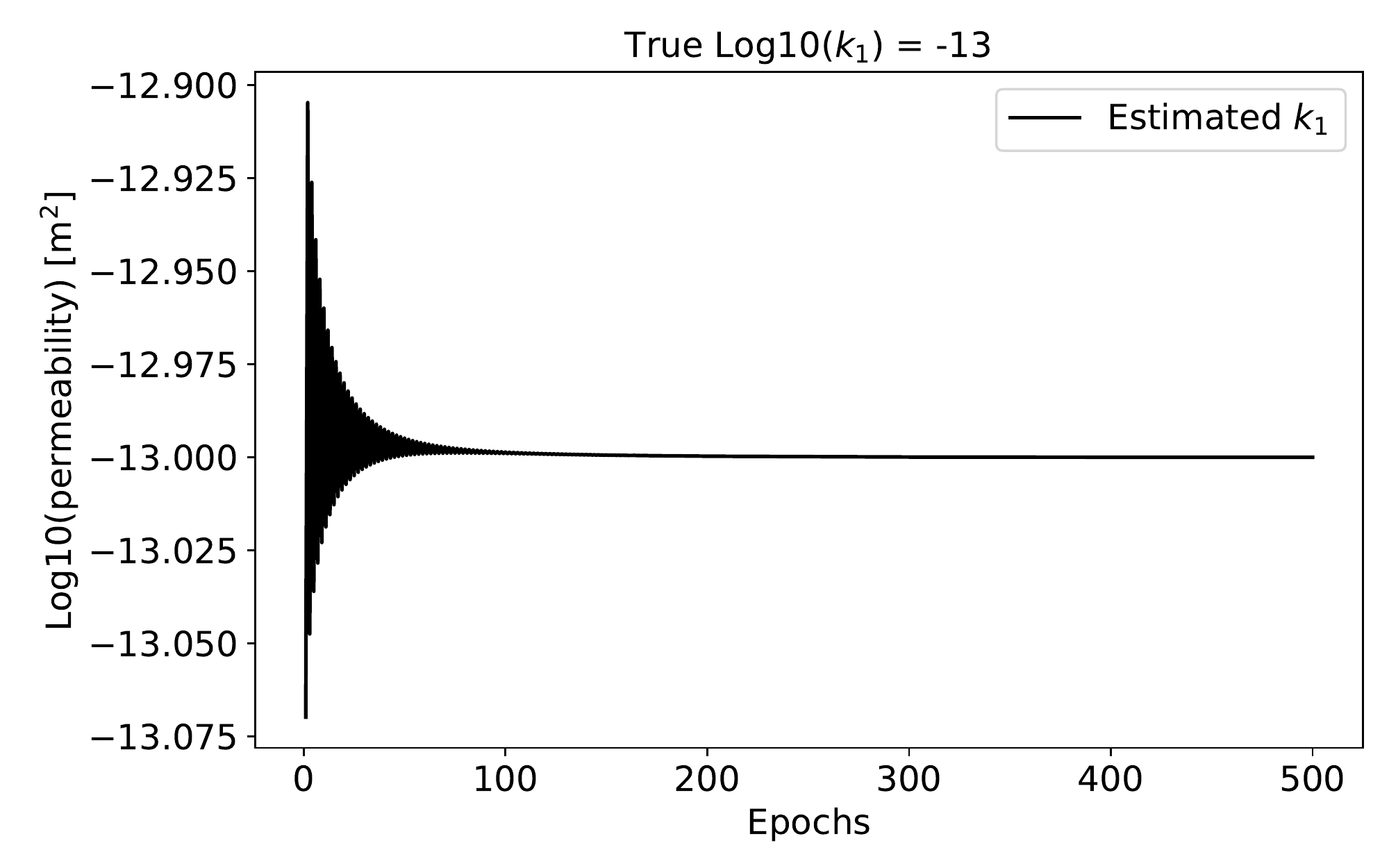}}
  \subfigure[$k_2$ training]
  {\includegraphics[width = 0.49\textwidth]
  {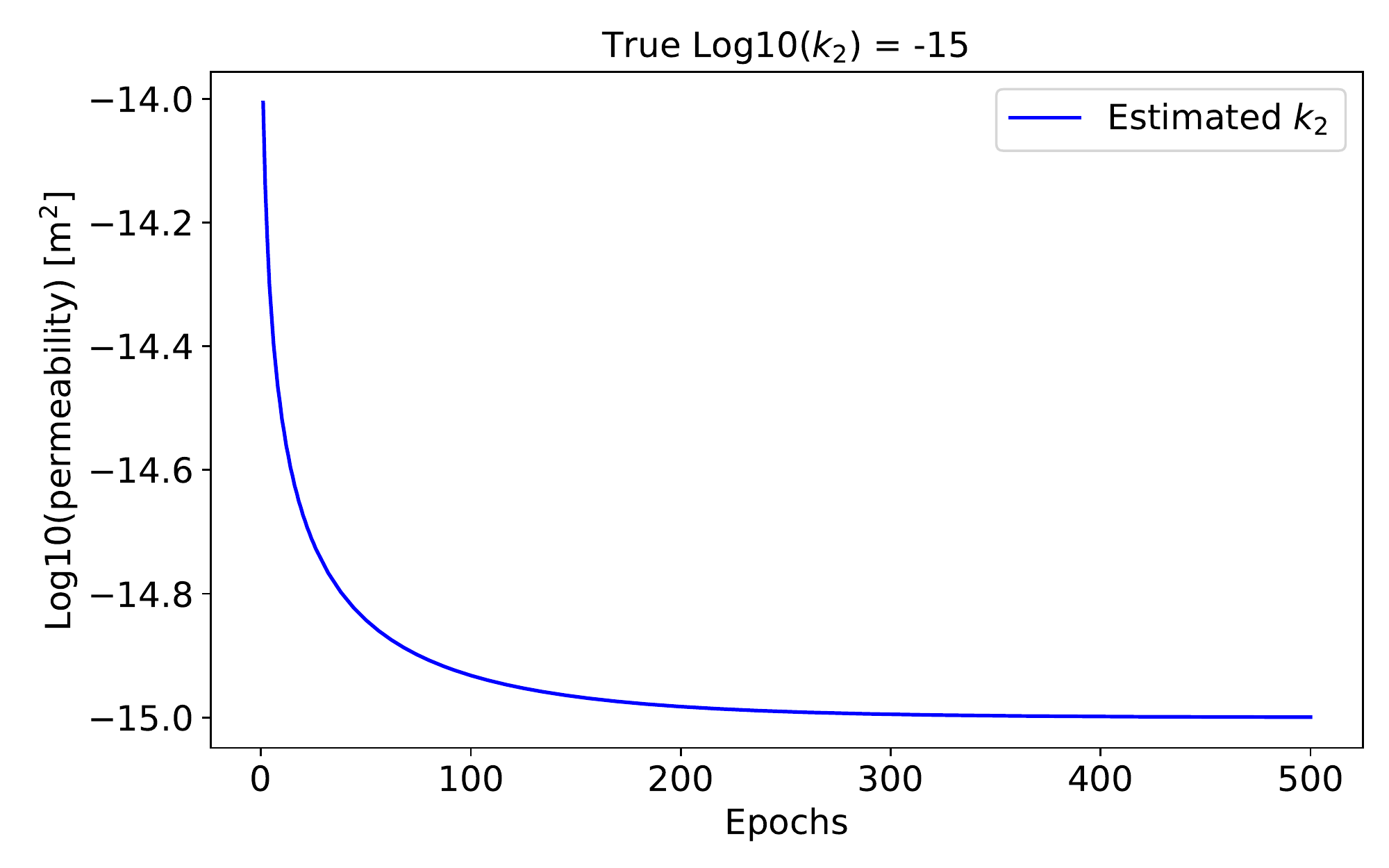}}

  \caption{\textbf{AdjointNet for inversion of heterogeneous permeability:}~This figure shows the training progress of the two permeability neural networks.
  We can see that the permeability estimates plateau within 100 epochs.} 
  \label{fig:hetero-prediction-vs-epochs}
\end{figure}

For this case study, the model domain consists of two zones with different permeabilities. 
The domain is meshed in the same manner as in Sec.~\eqref{subsec:homogeneous_case}.
The total number of grid cells is equal to 100.
The initial condition for the \texttt{PFLOTRAN} simulation includes the pressure of 1.5$\times10^5$ Pa throughout the model domain.
Pressure boundary conditions are prescribed on the left and right sides of the model domain, which are 1.0 and 0.5\, MPa.
For ground truth, we assume the first and last $50$ grid cells have permeability values of $10^{-13}$ and $10^{-15}$\,m$^2$, respectively.
Based on these permeability values, we generate pressure data for AdjointNet-enabled inversion. 
The solid red line in Fig.~\eqref{fig:hetero-prediction-vs-real} shows the ground truth pressure distribution. 
Two neural networks with two layers and 256 neurons in each layer are trained to estimate the permeabilities until a loss value of $10^{-8}$ is reached.
The training process converged within 100 epochs, as shown in Fig.~\eqref{fig:hetero-prediction-vs-epochs}.
The AdjointNet learned permeabilities were able to predict the pressure distribution accurately, as seen in Fig.~\eqref{fig:hetero-prediction-vs-real}.  

%

\subsection{Lid-driven cavity flow using the Navier-Stokes equation}
In this section, we show the applicability of AdjointNet methodology to simulate 2D lid-driven cavity flow.
The domain is $4 \times 4$\,\si{\square\metre} discretized by $41 \times 41$ uniform grid points.
We used a pressure Poisson finite difference solver written in Python \cite{Cavity_NS}. 
The initial condition is $u, v, p = 0$ everywhere.
The boundary conditions are: $u=1$ at $y=4$; $u, v=0$ at $y=0$; $\frac{\partial p}{\partial y}=0$ at $y=0$, $p=0$ at $y=4$, $\frac{\partial p}{\partial x}=0$ at $x=0,4$.
The ground truth is simulated using $\rho$ and $\nu$ to be equal to 1\,\si{kg/\cubic\metre} and 0.1\,\si{\square\metre/s}, respectively. 

\begin{figure}
  \centering
    {\includegraphics[width = 0.6\textwidth]
    {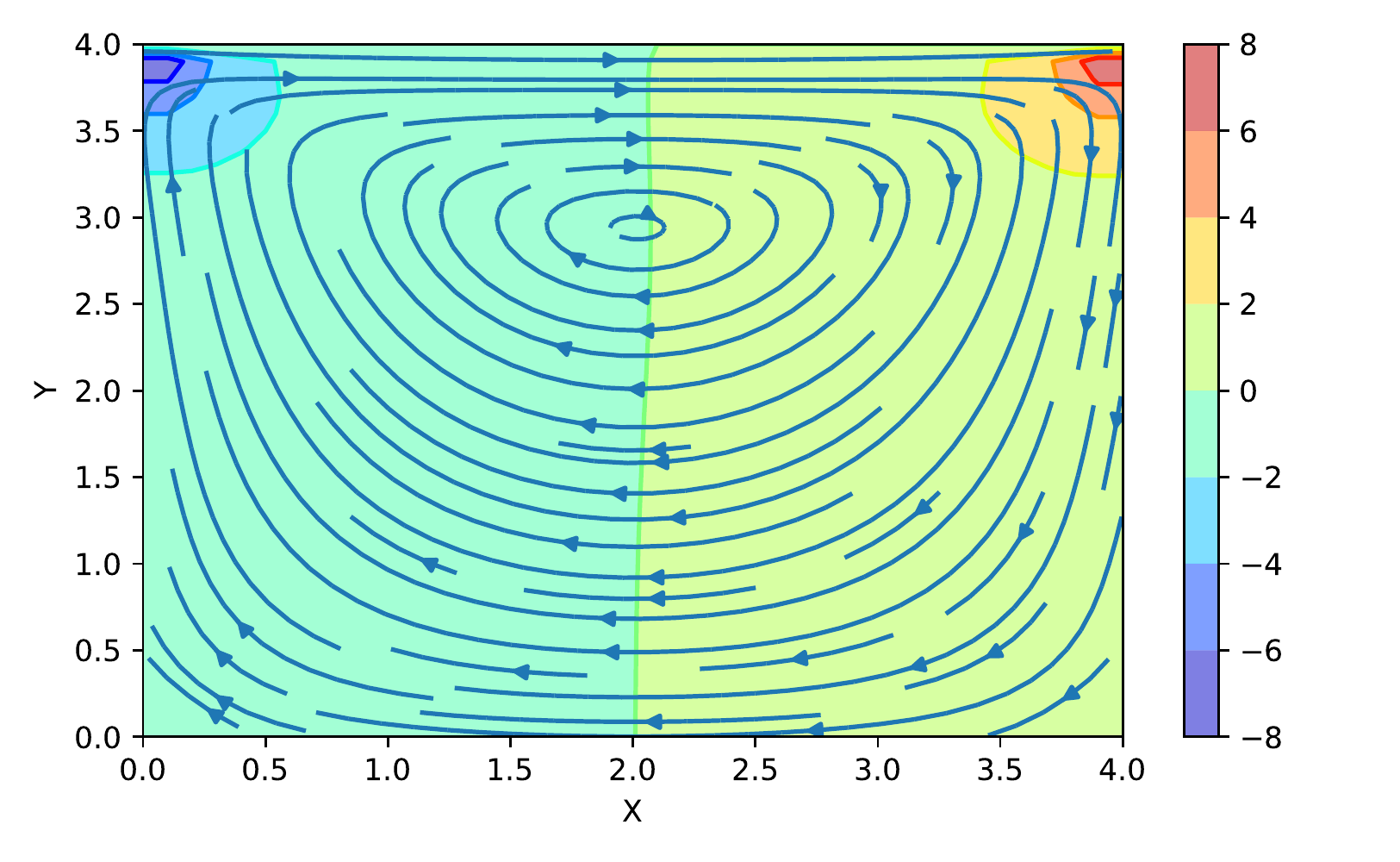}}
  \caption{\textbf{AdjointNet for inversion of kinematic viscosity:}~Lid-driven cavity Navier-Stokes flow. 
  Colors represent pressure distribution for the ground truth $\rho$ and $\nu$.}
  \label{fig:NS1}
\end{figure}

\begin{figure}
  \centering
    \subfigure[Loss versus epochs]
    {\includegraphics[width = 0.45\textwidth]
    {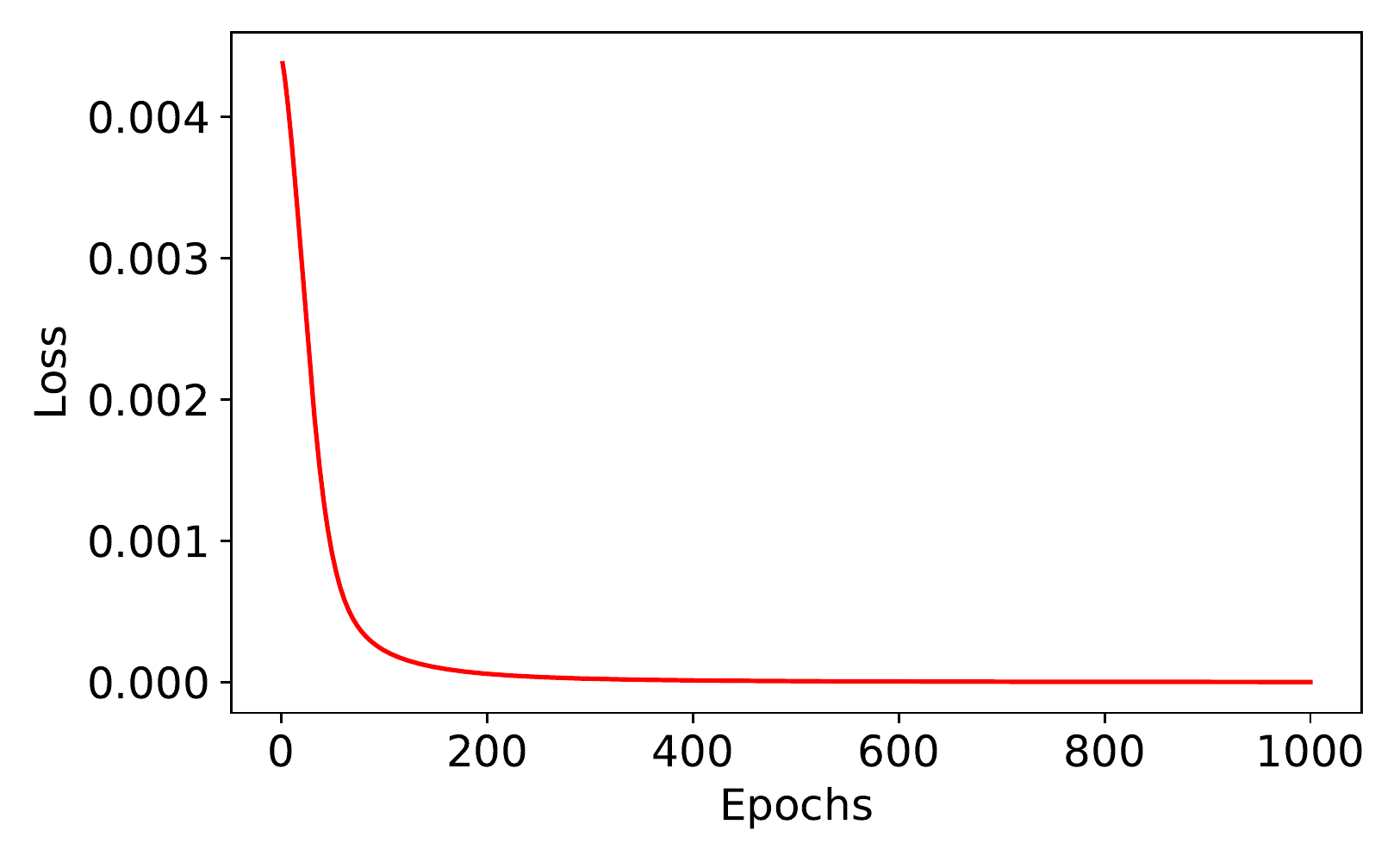}}
  \subfigure[True and estimated kinematic viscosity]
  {\includegraphics[width = 0.45\textwidth]
  {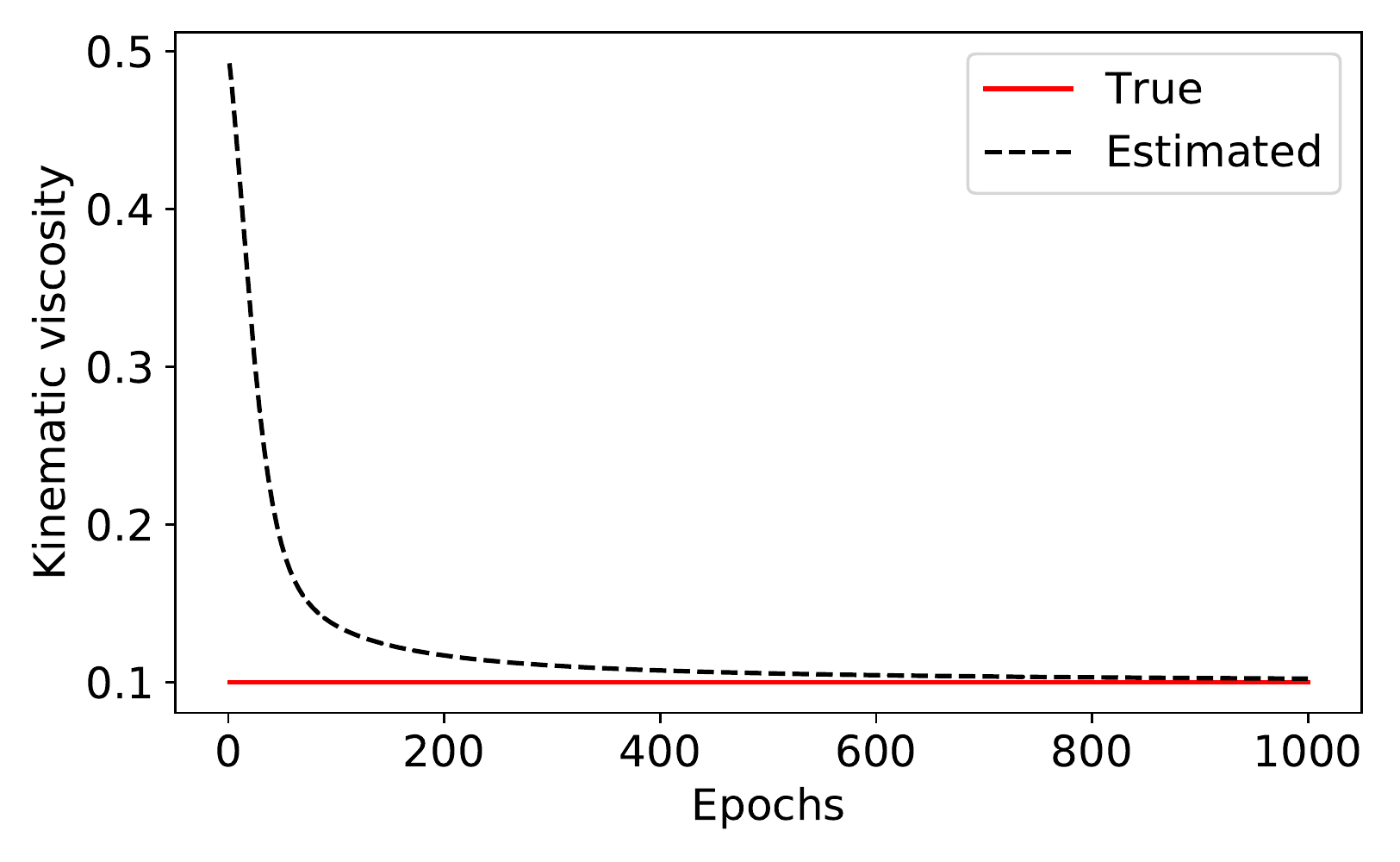}}

  \caption{\textbf{AdjointNet for inversion of kinematic viscosity:}~Kinematic viscosity prediction by AdjointNet (a) shows the progression of loss over epochs and 
  (b) compares the true and kinematic viscosity as training progresses. 
  \label{fig:NS2}}
\end{figure}

The model is run for 100 iterations, and for each iteration total time step is 300. 
Simulations are performed to generate pressure data for AdjointNet-enabled inversion.
A pressure snapshot (Figure~\eqref{fig:NS1}) with 1,681 pressure points is used to learn the kinematic viscosity.
A neural network for the kinematic viscosity is created that has two layers containing 64 and 32 neurons in each layer. 
Convergence is reached in around 200 epochs with the loss value of less than $10^{-8}$ (Fig.~\eqref{fig:NS2}).
AdjointNet estimated kinematic viscosity is 0.1021 compared to the ground truth of kinematic viscosity 0.1, which shows the applicability of AdjointNet for inversion of flow parameters with the Navier-Stokes equation.

\section{Conclusions}\label{sec:conclusions}
Currently, neural network models are black-box models, making it challenging to interpret and trust the results.
Recently, XAI techniques (e.g., deep Taylor decomposition) were developed to make neural networks interpretable, but they are not universally applicable (i.e., trustworthiness is always problem-specific).
This paper proposed a new ML methodology called AdjointNet and demonstrated its utility to solve inverse problems.
The proposed AdjointNet framework overcomes these pitfalls, as the existing codes are interpretable and physics is constrained everywhere.
Through representative numerical examples, we showed that one could invert material parameters by embedding a physics-code (e.g., \texttt{PFLOTRAN}, Navier-Stokes solver) in AdjointNet without modifying the underlying physics-code, or without the need to re-write a  code for a PDE from scratch, which is a significant bottleneck for existing PIML methods such as PINNs.
We demonstrated AdjointNet's capability on four examples: (1) homogeneous permeability inversion; (2) data assimilation; (3) heterogeneous permeability inversion; and (4) estimating material properties in a Navier-Stokes lid-driven cavity flow. 

In addition to satisfying the physics constraint, the results also showed that we could successfully interpret and trust them as the existing codes are verified and validated by domain experts over decades.
Since AdjointNet requires sensitivities of the underlying variable (e.g., pressure, temperature) with respect to the input parameters, this can be a limitation in cases where there are many parameters, such as spatially distributed permeability, in the case of porous flow. 
In such cases, numerically evaluating the sensitivities by perturbing the input parameters may be computationally intractable. One will have to resort to the discrete adjoint method, which is agnostic to the number of model parameters but may need some intervention to the underlying physics code. 
Alternatively, one can resort to \textit{a priori} dimensionality reduction, where one can integrate local/global sensitivity analysis with unsupervised learning (e.g., $k$-mean clustering) to identify the most critical parameters and estimate them using AdjointNet.
Finally, since AdjointNet can learn the underlying parameters while constraining the physics, it can perform well beyond the training range.
As a result, this framework is attractive for domain scientists who use complex codes to simulate physical processes. 

\section*{ABBREVIATIONS}
\begin{itemize}
  \item HPC:~High Performance Computing
  \item LIME:~Local Interpretable Model-Agnostic Explanations
  \item LRP:~Layer-Wise Relevance Propagation
  \item KGML:~Knowledge-Guided Machine Learning 
  \item MCMC:~Markov Chain Monte Carlo
  \item ML:~Machine Learning
  \item NMFk:~Non-negative Matrix Factorization with custom $k$-means Clustering
  \item NTFk:~Non-negative Tensor Factorization with custom $k$-means Clustering
  \item PDE:~Partial Differential Equation
  \item \texttt{PFLOTRAN}:~An open source, state-of-the-art massively parallel subsurface flow and reactive transport code
  \item PIML:~Physics-informed Machine Learning 
  \item PINN:~Physics-informed Neural Network
  \item XAI:~eXplainable Artificial Intelligence
\end{itemize}
%
\section*{Acknowledgments}
The authors thank the U.S. Department of Energy's Biological and Environmental Research Program for support through the SciDAC4 program. SK and BA also thank the Center for Space and Earth Science and the Information Science \& Technology Institute at Los Alamos National Laboratory. 
MKM thanks the support from the U.S. Department of Energy Office of Science's River Corridor Science Focus Area at PNNL.
The views and opinions of authors expressed herein do not necessarily state or reflect those of the United States Government or any agency thereof.

\newpage
\bibliographystyle{IEEEtran}
\bibliography{references}
\end{document}